# On a class of harmonic-like numbers

___________________________

**Definition**

The classical harmonic numbers (see for instance [1] for a definition and their main properties) $H_n = 1 + 1/2 + 1/3 + ..... + 1/n$ possess the integral representation

$$H_n = \int_0^1 \frac{1 - t^n}{1 - t} \, dt$$

from which a "natural" generalisation seems to be, for any complex number $a$ :

$$H_n(a) = \int_0^1 \frac{a^n - t^n}{a - t} \, dt \qquad (1)$$

In effect $a^n - t^n = (a - t)(a^{n-1} + a^{n-2} t + ..... + a^{n-p} t^{p-1} + .... t^{n-1})$, and, through a term by term integration, we obtain

$$H_n(a) = a^{n-1} + \frac{a^{n-2}}{2} + ........ + \frac{a^{n-p}}{p} + ...... + \frac{1}{n}$$

Then if we set $a = 1$ in this relation we find the usual harmonic numbers. Hence the label "harmonic-like" given to the numbers $H_n(a)$.

**An application to analysis**

We defer the study of the general number $H_n(a)$ to a further note, and instead we focus on the particular case $H_n(\tfrac{1}{2})$, to find an application to the expansion into series of some special functions.

To this end we need to transform the integral representation *(1)*.

By some obvious substitutions of variable in *(1)* we get:

$$H_n(\tfrac{1}{2}) = \int_0^1 \frac{\left(\tfrac{1}{2}\right)^n - t^n}{\left(\tfrac{1}{2}\right) - t} \, dt = \int_{-\tfrac{1}{2}}^{\tfrac{1}{2}} \frac{\left(\tfrac{1}{2}\right)^n - \left(\tfrac{1}{2} - y\right)^n}{y} \, dy = \frac{1}{2^n} \int_{-1}^1 \frac{1 - (1 - z)^n}{z} \, dz$$

and finally

$$H_n(\tfrac{1}{2}) = \frac{1}{2^n} \int_0^1 \frac{(1+x)^n - (1-x)^n}{x} dx \qquad (2)$$

Now we consider the integral sine function, classically defined [2] as

$$Si(z) = \int_0^z \frac{\sin x}{x} dx \qquad .$$

or equivalently

$$Si(z) = \int_0^1 \frac{\sin zx}{x} dx \qquad .$$

If we multiply both sides of this last relation by the first derivative of $Si(z)$, i.e., $(\sin z)/z$, then

$$\frac{d}{dz}\{Si(z)\}^2 = 2\int_0^1 \frac{\sin z \sin zx}{zx} dx \qquad .$$

But since $2\sin u.\sin v = \cos(u-v) - \cos(u+v)$, this leads to

$$\frac{d}{dz}\{Si(z)\}^2 = \int_0^1 \frac{\cos z(1-x) - \cos z(1+x)}{zx} dx \qquad .$$

Now by using the MacLaurin expansion of the $\cos u$ functions under the integral sign (with $u = z(1 - x)$ and $u = z(1 + x)$ successively), then the expression *(2)* of the numbers $H_n(\tfrac{1}{2})$, and integrating the result with regard to $z$, we finally obtain a power series expansion of the square of the integral sine function, valid for any complex number $z$,

$$\{Si(z)\}^2 = \sum_{n=1}^{+\infty} \frac{(-1)^{n+1}}{(2n)!(2n)} H_{2n}(\tfrac{1}{2}) (2z)^{2n} \qquad .$$

which contains $H_n(\tfrac{1}{2})$ coefficients of only even subscripts. On the other hand, by the same technique as above, we obtain

$$2\cos z\, Si(z) = \sum_{n=1}^{+\infty} \frac{(-1)^{n+1}}{(2n-1)!} H_{2n-1}(\tfrac{1}{2}) (2z)^{2n-1} \qquad .$$

which contains $H_n(\tfrac{1}{2})$ coefficients of only odd subscripts.

If we replace the trigonometric functions *sint* and *cost* by their hyperbolic counterparts *sht* and *cht* , the same techniques as above lead us to similar series expansions.

**Other property of $H_n(\tfrac{1}{2})$**

In passing, we notice that the numbers $H_n(\tfrac{1}{2})$ display a remarkable relationship with the **inverses of binomial coefficients**, for we have the following identity, several times independently proved ([3], [4] and [5]), and even generalised ([6])

$$\sum_{p=1}^{n} \binom{n}{p}^{-1} = \frac{n+1}{2^{n+1}} \sum_{k=1}^{n+1} \frac{2^k}{k} \quad .$$

According to our definition the right-hand side of this relation is equal to $(n+1)H_{n+1}(\tfrac{1}{2})$. Therefore, for $n \geq 1$ we have the simple identity:

$$H_{n+1}(\tfrac{1}{2}) = \frac{1}{n+1} \sum_{p=1}^{n} \binom{n}{p}^{-1} \quad .$$

________________________

*Juan PLA, 315 rue de Belleville 75019 Paris (France)*